\documentclass[12pt]{amsart}
\usepackage{amsmath, amssymb}

\RequirePackage[top=4.0cm,bottom=3.0cm,inner=3.0cm,outer=3.0cm,
headsep=1cm,footskip=2cm]{geometry}
\RequirePackage{ifpdf}
\RequirePackage{xr-hyper}
\RequirePackage[colorlinks=true]{hyperref}

\newtheorem{thm}{Theorem}

\newtheorem{prop}[thm]{Proposition}
\newtheorem{cor}[thm]{Corollary}
\newtheorem{defi}[thm]{Definition}
\newtheorem{rmk}[thm]{Remark}

\newcommand{\thmref}[1]{Theorem~\ref{#1}}

\begin{document}

\title
{Sign changes of coefficients of certain Dirichlet series}

\author[J. Meher, S. Pujahari and K. Deo Shankhadhar]{Jaban Meher, Sudhir Pujahari and Karam Deo Shankhadhar }

\address[Jaban Meher]{Department of Mathematics, Queen's University,
Kingston, Ontario, K7L 3N6, Canada\\
{\sl Current Address:} Department of Mathematics, Indian Institute of Science, Bangalore 560012, India}
\email{jaban@math.iisc.ernet.in}

\address[Sudhir Pujahari]{Department of Mathematics, Indian Institute of Science 
Education and Research, 900, NCL Innovation Park, Dr Homi Bhabha Road,
Pune 411008, India}
\email{sudhir.pujahari@students.iiserpune.ac.in}

\address[Karam Deo Shankhadhar]{The Institute of Mathematical Sciences,
IV cross road, CIT campus, Taramani,
Chennai 600113, India\\
{\sl Current Address:} Depto. Matematicas, Facultad de Ciencias, Univ. de Chile, Las Palmeras 3425, \~Nu\~noa, Santiago, Chile}
\email{karam@imsc.res.in, karamdeo@u.uchile.cl}




\maketitle

\begin{abstract}
In this paper, we give criteria for infinitely many sign changes of the coefficients of any Dirichlet series 
if the coefficients are real numbers. We also provide examples where our criteria are applicable.
\end{abstract}

\section{Introduction}

The sign changes in a sequence of real numbers have attracted many for different reasons.
In particular, the signs of the Fourier coefficients of different automorphic forms have been studied
due to their various number theoretic applications. Coming back to the general scenario, given a 
sequence of real numbers, one may ask a natural question about the number of sign changes in 
the sequence. Pribitkin \cite{pribitkin1} proved a general result about the sign changes of the 
coefficients of any general Dirichlet series. More precisely, he proved that the coefficients of 
a general Dirichlet series change signs infinitely often. In \cite{pribitkin2}, he has given numerous 
applications of his result. In this article, we prove two results about the sign changes of
the coefficients of certain kind of Dirichlet series. Using the method similar to Kohnen's method of proof in 
\cite{kohnen}, we prove
our first result about the sign changes. In our next result, we give a criterion for a Dirichlet series
with real coefficients to have infinitely many sign changes. To prove our second result, we follow the method 
similar to the proof 
of the main result of \cite{hulse}. Furthermore, we provide applications of these two results in \S \ref{applications}. 
Using these results,
we conclude the oscillatory behaviour (infinitely many sign changes) of the Fourier coefficients of elliptic cusp forms, 
Siegel cusp forms, 
Maass cusp forms and second order cusp forms. 
The sign change results in the case of elliptic cusp forms, Siegel cusp forms and symmetric power ~$L$-functions associated to elliptic cusp forms
have been established in \cite{KKP, K07, MSV} using different techniques. The applications provided in the case of second order cusp forms and 
Maass cusp forms are new. We establish the oscillatory behaviour of the Fourier coefficients of second order cusp forms and infinitely many sign change
in the subsequences ~$\{a_{n^j}\}_{n=1}^{\infty}$~ of the sequence ~$\{a_{n}\}_{n=1}^{\infty}$~of the Fourier coefficients of Maass cusp forms for any ~$j=1,2,3,4$.

J. Kaczorowski kindly informed the third author that his recent paper with A. Perelli  \cite{KP} deals with a similar subject as this though the methods used are quite different.
More precisely, under some simple conditions they show that the real and imaginary parts of any linear combination of coefficients of $L$-functions from the Selberg class have infinitely many sign changes by proving a general $\Omega$-theorem for the coefficients of polynomial combinations of such $L$-functions.

\section{Statement of results}

Let $k$ and $N$ be given positive integers. Put $W_N=\begin{pmatrix}0&-1\\N&0
\end{pmatrix}$. Let $f$ be any non-zero function defined on the upper half-plane 
$\mathcal{H}=\{z\in \mathbb{C} : \mbox{Im}(z)>0\}$. Define
$f|_kW_N(z):= (\sqrt{N}z)^{-k}f(\frac{-1}{Nz})$. Our first result is the following.
\begin{thm}\label{fourier-sign}
Suppose that $f$ and $g=f|_kW_N$ have the following 
Fourier series expansions.
$$f(z)=\sum_{n\ge 0}a_ne^{2\pi inz} ~{\rm{and}}~ g(z)=f|_kW_N(z)=\sum_{n \ge 0}b_ne^{2 \pi inz}.$$
Assume that the Fourier coefficients $a_n, b_n ~(n\ge 1)$ are bounded by $O(n^\alpha)$, 
where $\alpha$ is a positive constant. 
If $a_0$ (respectively $b_0$) is zero and $b_n$ (respectively $a_n$) are real numbers for $n\ge 1$, then the sequence  
${\{b_n\}}_{n=1}^{\infty}$ (respectively  
${\{a_n\}}_{n=1}^{\infty}$) changes signs infinitely often.
\end{thm}
As a consequence of the above theorem, we have the following corollary.
\begin{cor}
Let $f$ and $g$ be as in the above theorem. Suppose that $a_0=b_0=0$ and $a_n, b_n$ are real numbers for $n\ge 1$. 
Then both the sequences ${\{a_n\}}_{n=1}^{\infty}$ and  
${\{b_n\}}_{n=1}^{\infty}$ change signs infinitely often.
\end{cor}
\begin{rmk}
The type of function considered in \thmref{fourier-sign} has appeared in the recent work of 
Choie and Kohnen \cite{choie-kohnen1}.
\end{rmk} 
We 
state our second result about the sign changes in a given sequence of real 
numbers.
\begin{thm}\label{signchange}
Let ${\{a_n\}}_{n=1}^{\infty}$ be a sequence of real numbers such that 
$a_n=O(n^{\alpha})$, for some real 
number $\alpha$. Assume that the Dirichlet series 
$\sum_{n=1}^{\infty}\frac{a_n}{n^s}$ can be analytically
continued to Re$(s)>r\ge0$ and has polynomial growth in Im$(s)$ in this region.
Furthermore, assume that the Dirichlet series 
$\sum_{n=1}^{\infty}\frac{{a_n}^2}{n^s}$ has a singularity at $s=k>0$ 
such that ~$\alpha + r<k$. Then the 
sequence ${\{a_n\}}_{n=1}^{\infty}$ changes signs infinitely often. 
\end{thm}
\begin{rmk}
The condition in the hypothesis of Theorem \ref{signchange} that 
~$\alpha+ r<k$~ is essential. 
For example, take the case 
of Riemann zeta function $\zeta(s)=\sum_{n\ge 1}\frac{1}{n^s}$, 
we have $\alpha =0, r=1$ and $k=1$. The condition $\alpha+r<k$ is not satisfied 
and all the coefficients are of positive sign.
\end{rmk}

\section{Proofs}

We first begin with a theorem due to Hecke \cite[Theorem 7.3]{iwaniec}.
\begin{thm}\label{hecke}
Let $k$ and $N$ be positive integers.
Suppose that $f$ and $g$ are given by the
Fourier series
$$f(z)=\sum_{n\ge 0}a_ne^{2\pi inz} ~{\rm{and}}~ g(z)=\sum_{n \ge 0}b_ne^{2 \pi inz},$$
with coefficients $a_n, b_n ~(n\ge 1)$ bounded by $O(n^\alpha)$, where $\alpha$ is a positive constant.
Let 
$$ L(s,f)=\sum_{n\ge 1}\frac{a_n}{n^s},~ L(s,g)=\sum_{n\ge 1}\frac{b_n}{n^s}$$
and
$$L^*(s,f)={(\sqrt{N}/{2 \pi})}^s\Gamma(s)L(s,f),~ L^*(s,g)={(\sqrt{N}/{2 \pi})}^s\Gamma(s)L(s,g).$$
Then the following assertions are equivalent:
\begin{enumerate}
\item[(i)]
The functions $f$ and $g$ are related by $g=f|_kW_N$, where $f|_kW_n$ is defined
just before the statement of Theorem \ref{fourier-sign}.
\item[(ii)]
Both $L^*(s,f)$ and $L^*(s,g)$ have meromorphic continuation over the whole complex plane,
$$
L^*(s,f)+\frac{a_0}{s}+i^{-k}\frac{b_0}{k-s},~~~~~L^*(s,g)+\frac{b_0}{s}+i^{-k}\frac{a_0}{k-s}
$$
are entire and bounded on vertical strips, and they satisfy
$$
L^*(s,f)=i^kL^*(k-s,g).
$$
\end{enumerate} 
\end{thm}
\vspace{.2cm}
\proof[Proof of Theorem \ref{fourier-sign}]~
Let ~$ L(s,f)=\sum_{n\ge 1}\frac{a_n}{n^s}$ ~and~ $L(s,g)=\sum_{n\ge 1}\frac{b_n}{n^s}$~
be the Dirichlet series associated to $f$ and $g$ respectively. Let
$$L^*(s,f)={(\sqrt{N}/{2 \pi})}^s\Gamma(s)L(s,f) ~~~\mbox{and}~~~ L^*(s,g)={(\sqrt{N}/{2 \pi})}^s\Gamma(s)L(s,g).$$
Since $g|_kW_N = (f|_kW_N)|_kW_N = (-1)^kf$, it is sufficient to consider the 
case $a_0 =0$ and $b_n(n\ge 1)$ are real. Let $a_0=0$. Assume on the contrary that ~$b_n\ge 0$~ for all but 
finitely many $n$. Then by Landau's
theorem on Dirichlet series with non-negative coefficients, $L(s,g)$ is either convergent for all 
$s\in \mathbb{C}$ or it has a singularity at the real point of the abscissa of convergence. 
Since $g= f|_kW_N$, using \thmref{hecke}, $L^*(s,g)+\frac{b_0}{s}$ is an entire function of ~$s$ and therefore 
$L^*(s,g)$ does not have any singularity in the region
Re$(s)>0$. Since $\Gamma(s)$ has its poles exactly 
at non-positive integers, $L(s,g)$ is an entire function of ~$s$. From the above discussion, the series $L(s,g)$ 
converges for all $s$ 
and it has zero at all negative integers ~$\{-1,-2,\ldots\}$. Let ~$m_1,m_2,\ldots, m_t$~ be all positive integers with
~$m_1<m_2 \ldots <m_t$~ such that ~$b_{m_1},b_{m_2},\ldots,b_{m_t}$~ are strictly negative. For
$\nu =-1,-2, \ldots$, we have
$$
\sum_{n\ge 1}\frac{b_n}{n^\nu}=0.
$$
Dividing both sides by $m_t^{-\nu}$, we get
\begin{equation}\label{mid}
\sum_{n\ge 1, \atop n\ne m_1, m_2, \dots, m_t} b_n{\left(\frac{n}{m_t}\right)}^{-\nu}
=-b_{m_1}{\left(\frac{m_1}{m_t}\right)}^{-\nu} - \dots -b_{m_t}.
\end{equation}
Letting $-\nu \rightarrow \infty$, we see that the right hand side of the above equation has the limit $-b_{m_t}$, 
which is strictly positive. If there were $n>m_t$ for which $b_n>0$, then the left hand side of \eqref{mid}
would tend to infinity as $-\nu \rightarrow \infty$, giving a contradiction. Thus, we deduce that ~$b_n=0$~ 
for all ~$n>m_t$.
Therefore for all ~$\nu=-1,-2,\ldots$,~ we have
\begin{equation}\label{1}
\sum_{n=1}^{m_t}\frac{b_n}{n^\nu}=0.
\end{equation}
Let $1\le m\le m_t$ be the largest positive integer for which $b_m\ne 0$. Then from \eqref{1}, we have
$$
\sum_{n=1}^{m-1}b_n{\left( \frac{n}{m}\right)}^{-\nu} + b_m =0.
$$
Again letting $-\nu \rightarrow \infty$, we get $b_m=0$. This is a contradiction to the assumption that $m$ is the
largest integer for which $b_m\ne 0$. This proves the theorem.
\qed 
\begin{rmk}
The above theorem can be proved also by using Theorem \ref{hecke} and
the result on 
Dirichlet series proved in \cite{pribitkin1}.
\end{rmk}
\proof[Proof of Theorem \ref{signchange}]
Assume on the contrary that $a_n\ge0$ for all $n>T$ for some sufficiently large
number $T$.
Since $a_n=O(n^{\alpha})$, the series $M(s)=\sum_{n=1}^{\infty}\frac{a_n}{n^s}$ is absolutely convergent
for $Re(s)>\alpha+1$. Using the inverse Mellin transform, we get
\begin{equation}\label{begin}
\frac{1}{2\pi i}  \int_{\alpha+1-i\infty}^{\alpha+1+i\infty} \! M(s) \Gamma(s)x^s  \mathrm{d} s
=\sum_{n\ge 1}a_ne^{-n/x}. 
\end{equation}
Since $M(s)$ can be analytically continued to $Re(s)>r\geq 0$ (which is assumed), for any 
$\epsilon>0$ if we move the line of integration to $Re(s)=r+\epsilon$, we see that on the vertical line
$Re(s)=r+\epsilon$, the gamma function decreases exponentially and the function $M(s)$ has polynomial 
growth.  Thus, we do not get any singularity for $Re(s)>r$. Then the integral on the left hand side of \eqref{begin}
is $O(x^{r+\epsilon})$, for any $\epsilon>0$. This implies that
\begin{equation}\label{mellin}
\sum_{n\ge 1} a_ne^{-n/x}=O(x^{r+\epsilon}).
\end{equation}
We have
\begin{equation*} 
e|\sum_{n>T}a_ne^{-n/x}|-e|\sum_{n\le T}a_ne^{-n/x}|\le e|\sum_{n\ge 1}a_ne^{-n/x}|\le \beta x^{r+\epsilon},
\end{equation*}
where $\beta$ is some constant and $e$ denotes the exponential. 
Since $a_n$ is non-negative for $n>T$, the above equation gives us the 
following.
\begin{equation}\label{final_1}
\sum_{T<n\leq x}a_n\le e\sum_{n>T}a_ne^{-n/x}\le \gamma x^{r+\epsilon},
\end{equation}
for some constant $\gamma$.

\noindent Now, we claim that if $A(x)=\displaystyle\sum_{T<n\le x}a_n^2$, then for every ~$c$~ with $0<c<k$, every constant
$\alpha_1>0$ and every $x$, there exists an $x_0>x$ such that
\begin{equation}\label{bigger}
A(x_0)\ge \alpha_1x_0^c.
\end{equation}
Assume towards a contradiction that there exists $c<k$ such that $A(x)=O(x^c)$. Using the partial summation 
formula, we have for $Re(s)>k$,
\begin{equation}\label{final_2}
\sum_{n>T}\frac{a_n^2}{n^s}=s \! \int_T^{\infty} \frac{A(u)}{u^{s+1}}\ \mathrm{d} u.
\end{equation}
Since $A(x)=O(x^c)$, the right hand side of \eqref{final_2} is an analytic function for $Re(s)>c$. Under the 
hypothesis of the theorem, the
left hand side of \eqref{final_2} has a singularity at $s=k$, which is a contradiction. This proves the claim 
\eqref{bigger}.

\noindent Using \eqref{bigger} and \eqref{final_1} with the assumption that $a_n=O(n^\alpha)$, we get the following. 
\begin{equation}
\alpha_1x_0^{-\alpha+c}\le x_0^{-\alpha}\sum_{T<n\le x_0}a_n^2\le \sum_{T<n\le x_0}\frac{a_n^2}{n^{\alpha}}
\le \lambda \sum_{T<n\le x_0} a_n\le \lambda \gamma x_0^{r+\epsilon}, 
\end{equation}
for some constant $\lambda$. This implies,
\begin{equation}
x_0^{c-r-\alpha-\epsilon}\le \frac{\lambda \gamma}{\alpha_1}.
\end{equation}
Since $r+\alpha<k$ holds by assumption, we can choose $c$ and $\epsilon$ in such a way that the exponent on the left hand side of the 
above equation is greater than $0$, whereas the right hand side would be less than $1$, giving a contradiction. 
This proves the theorem.
\qed

\section{Applications}\label{applications}

\subsection{Applications of Theorem \ref{fourier-sign}}
\begin{enumerate}
\item[(i)]
{\bf{Elliptic cusp forms}}\\
Let ~$k, N$~ be positive integers. Let us denote the space of cusp forms (elliptic) of weight $k$ and 
level $N$
by $S_k(N)$. If $f$ is such a cusp form then it is easy to see that $f|_k W_N \in S_k(N)$. 
Therefore the Fourier coefficients of both $f$  and $f|_k W_N$ are bounded by the Hecke's trivial bound
$O(n^{k/2})$. 
Assume that the Fourier coefficients of $f$ are
real numbers, then applying
\thmref{fourier-sign}, we deduce 
the following proposition.
\begin{prop}\label{prop:1.1}
If $f$ is any elliptic cusp form of weight $k$ and level $N$ 
with real Fourier coefficients, then the coefficients change signs infinitely often.
\end{prop}
\item[(ii)]
{\bf{Second order cusp forms}}
\begin{defi}\label{def:1.2}
A holomorphic function $f:\mathcal{H}\longrightarrow \mathbb{C}$ is called a second order cusp form of weight $k$ and 
level $N$ if it satisfies the following conditions.
\begin{enumerate}
\item[(a)] The function $f|_k(\gamma-1)(z):=(cz+d)^{-k}f\left(\frac{az+b}{cz+d}\right)-f(z) \in S_k(N)$ for all 
$\gamma=\begin{pmatrix}a&b\\c&d\end{pmatrix} \in \Gamma_0(N)$.
\item[(b)] $f|_k(\pi-1)=0$ for all parabolic elements $\pi$ in $\Gamma_0(N)$.
\item[(c)] $f$ has exponential decay at each cusp of $\Gamma_0(N)$.
\end{enumerate}
\end{defi}
\noindent We denote the space of second order cusp forms of weight $k$ and level $N$
by $S_k^2(N)$. It is proved in \cite[Proposition 10]{secondorder} that
if $f\in S_k^2(N)$ then $f|_kW_N \in S_k^2(N)$.  
The $n$-th Fourier coefficient of any second order cusp form of weight $k$ and level $N$ satisfies the 
trivial bound $O(n^{k/2}\log n)$ \cite[Lemma 7]{secondorder}. We get the following result by applying
\thmref{fourier-sign} in this situation. 
\begin{prop}
If the Fourier coefficients of a second order cusp form are real numbers, then they change signs infinitely often.
\end{prop}
\end{enumerate}

\subsection{Applications of Theorem \ref{signchange}}
\begin{enumerate}
\item[(i)]
{\bf{Elliptic cusp forms}}\\
Let ~$f(z)=\sum_{n\ge 1}a(n)e^{2\pi inz}$~ be any cusp form of weight $k$ and level
$N$ such that $a(n)$ are 
real numbers. We normalize the
Fourier coefficients by setting $a_n=a(n)/n^{(k-1)/2}$. We see that, the sequence of real numbers 
~$\{a_n\}_{n=1}^{\infty}$~satisfies the hypothesis of \thmref{signchange}. We have the Deligne's 
bound ~$a_n=O(n^{\epsilon})$, for any $\epsilon >0$. The series  
$\sum_{n=1}^{\infty}\frac{a_n}{n^s}$ extends to an 
analytic function to the whole complex plane and has a polynomial growth in Im$(s)$  
in the region Re$(s)>0$. Using Rankin-Selberg
convolution, one can prove that the series ~$\sum_{n=1}^{\infty}\frac{{a_n}^2}{n^s}$~ has a pole at $s=1$. 
In this case, ~$\alpha +r =\epsilon < 1$.
Applying 
\thmref{signchange}, we get infinitely many sign changes for the sequence 
${\{a_n\}}_{n=1}^{\infty}$. Note that the Hecke's trivial bound of ~$a(n)=O(n^{k/2})$~ would also suffice to use \thmref{signchange} to get this application. 
\item[(ii)]
{\bf{Siegel cusp forms of degree $2$}}\\
Let $F$ be a degree $2$ Siegel cusp form of weight $k$ 
for the Siegel modular group
~$Sp_4(\mathbb{Z})$~, which is an eigenform for all the Hecke operators ~$T(n)$~
with eigenvalues ~$\lambda_F(n)$~ 
and not a Saito-Kurokawa lift. Then it is well known that $\lambda_F(n)$ are real numbers. We normalize the
eigenvalues by setting $\lambda_n=\lambda_F(n)/n^{k-3/2}$. We know that $\lambda_n=O(n^\epsilon)$ for 
any $\epsilon >0$, as a consequence of the Ramanujan-Petersson conjecture, proved by Weissauer \cite{weissauer}. 
The analytic properties of the Dirichlet series $\sum_{n\ge 1}\lambda_nn^{-s}$ are well studied. In fact, it can be 
analytically continued to Re$(s) >0$ and has polynomial growth in Im$(s)$ in this region (see \cite{kohnen-sengupta}). 
In \cite[Theorem 1.1]{DKS},
it has been proved that the Dirichlet series $\sum_{n\ge 1}\lambda_n^2n^{-s}$ admits an analytic continuation 
to Re$(s)>1/2$ with the exception of a simple pole at $s=1$. In this case,
~$\alpha+r = \epsilon <1$~. Thus the sequence of eigenvalues
$\{\lambda_n\}_{n=1}^{\infty}$ satisfies all the hypothesis of Theorem \ref{signchange} and therefore it is 
oscillatory (infinitely many sign changes). 
\item[(iii)]
{\bf{Maass cusp forms}}\\
Let $g(z)$ be a Maass cusp form for the full modular group $SL_2(\mathbb{Z})$
with Laplace eigenvalue $1/4+\nu^2$. Suppose that the 
Fourier expansion of ~$g$~ is given by
 $$
 g(z)=\sum_{n\ne 0}a_n\sqrt{y}K_{i\nu}(2\pi |n|y)e^{2\pi inx},
 $$
where $z=x+iy$ and $K_{i\nu}$ is the modified Bessel function of the third kind.
Assume that $a_n~(n\ge 1)$ are real numbers.
The Dirichlet series attached to $g$ is defined as follows.
$$
L(s,g)=\sum_{n\ge 1}\frac{a_n}{n^s}.
$$
By a result of Kim and Sarnak \cite{kim-sarnak}, 
we have ~$a_n=O(n^{7/64+\epsilon})$, for any $\epsilon >0$. Let us fix ~$\mu=0$~ or ~$1$~
accordingly $g$ is even or odd respectively. Let
$$
\Lambda(s,g)=\pi^{-s}\Gamma\left(\frac{s+\mu+i\nu}{2}\right)\Gamma\left(\frac{s+\mu-i\nu}{2}\right)L(s,g).
$$
It is well known that $\Lambda(s,g)$ is entire and hence $L(s,g)$ is entire. The Dirichlet series $L(s,g)$ 
has polynomial growth in Im$(s)$ for Re$(s)>0$. By Rankin-Selberg convolution, the series 
$\sum_{n\ge 1}\frac{{a_n}^2}{n^s}$ has a simple pole at $s=1$. In this case,
~$\alpha+r = 7/64+\epsilon < 1$. Applying \thmref{signchange}, we deduce that
the sequence ${\{a_n\}}_{n=1}^{\infty}$ changes signs infinitely often.
\begin{rmk}
The sign change result for elliptic cusp forms is obtained in \cite{KKP} also.
The result pertaining to Siegel cusp forms of degree ~$2$~ was first established by Kohnen in 
\cite{K07}. In the case of Maass cusp forms, we do not know whether the Dirichlet series
attached to it, 
has infinitely many real zeros or not if $\nu\ne 0$. Thus, we cannot apply Theorem $1$ of \cite{pribitkin1} in this case.
However, in \cite{KKP} it is remarked that
the sign change result for Maass cusp forms can be formulated along the lines of their theorems for elliptic cusp forms.
\end{rmk}
\item[(v)] 
{\bf Symmetric power $L$-functions associated to elliptic cusp forms}\\
Let us assume that ~$f(z)=\sum_{n\ge 1}a(n)e^{2\pi inz} \in S_k(1)$~ be a Hecke eigen cusp form (that is, an eigenfunction for
all the Hecke operators ~$T(n), n\geq 1$) of weight $k$ for the
full modular group $SL_2(\mathbb{Z})$. Let $a_n=\frac{a(n)}{n^{(k-1)/2}}$. It is known that ~$a_n, n\geq 1$~
are real numbers. For any ~$j=2,3,4$, using \thmref{signchange}
we show that the subsequence {\bf{~$\{a_{n^j}\}_{n=1}^{\infty}$~}} changes signs infinitely often. For any fix 
~$j\in\{2,3,4\}$, we verify that the sequence of real numbers ~$\{a_{n^j}\}_{n=1}^{\infty}$~
satisfies the hypothesis of \thmref{signchange}. For any ~$n\geq 1$ we have ~$a_{n^j}=O(n^{j\epsilon})$~ for any
~$\epsilon >0$. Consider the Dirichlet series ~$\sum_{n\geq 1}a_{n^j} n^{-s}$. For ~$j=2$, this series can be analytically
continued to the half plane Re$(s) > \frac{1}{2}$~ \cite[\S 13.8]{iwaniec} and has polynomial growth in this region. 
From \cite[(2.7),(2.10), Proof of Theorem 1.2]{LU}, we have
the analytic continuation and polynomial growth in Im$(s)$ of the series ~$\sum_{n\geq 1}a_{n^j} n^{-s}, j=2,3$~ in the region 
Re$(s) > \frac{1}{2}$. From \cite[(3.1)]{LS}, we know that the series ~$\sum_{n\geq 1}a_{n^j}^2 n^{-s}$~ has a simple pole
at ~$s=1$. For any ~$j=2,3,4$, ~$\alpha+r=j\epsilon+\frac{1}{2} < 1$. Applying \thmref{signchange} we get the infinitely many
sign changes for any subsequence ~$\{a_{n^j}\}_{n=1}^{\infty}, j=2,3,4$. 
\begin{rmk}
In \cite{MSV}, certain quantitative results for the sign change in each subsequence ~$\{a_{n^j}\}_{n=1}^{\infty}, 
j\in\{2,3,4\}$~ have been established and hence infinitely many sign change in the subsequences. The proof uses suitable bounds for the two average sums
~$\sum_{1\leq n\leq x}a(n^j)$~ and ~$\sum_{1\leq n \leq x}{a(n^j)}^2$. Here we are using less tools to
conclude the infinitely many sign change but we do not get any quantitative
result.
\end{rmk}
\item[(vi)]
{\bf Symmetric power $L$-functions associated to Maass cusp forms}\\
Let ~$g(z)= \sum_{n\ne 0}a_n\sqrt{y}K_{i\nu}(2\pi |n|y)e^{2\pi inx}$~ be a Maass cusp form for the full modular group $SL_2(\mathbb{Z})$
with Laplace eigenvalue $1/4+\nu^2$. Suppose that it is an eigenfunction for all the Hecke operators ~$T(n), n\geq 1$. 
Using \thmref{signchange}, we show that the subsequences {\bf{~$\{a_{n^j}\}_{n=1}^{\infty}$~}}, $j=2,3,4$ change signs infinitely often.
For any ~$n\geq 1$ we have ~$a_{n^j}=O(n^{j(\frac{7}{64}+\epsilon)})$~ for any ~$\epsilon >0$.
Consider the Dirichlet series ~$L_2(s):=\sum_{n\geq 1}a_{n^2} n^{-s}$. Let ~$L(s,{\rm sym}^2g)$~ denotes the symmetric square
$L$-function attached to ~$g$. It is known that, 
$$
L(s,{\rm sym}^2g) = \zeta(2s) \sum_{n\geq 1}a_{n^2} n^{-s},  ~~~~~~~~{\rm Re}(s) >1.
$$
Since ~$L(s,{\rm sym}^2g)$~ can be extended to an entire function, the series ~$\sum_{n\geq 1}a_{n^2} n^{-s}$~
admits an analytic continuation for the region Re$(s) >1/2$. 
For ~$j=3,4$, consider the Dirichlet series ~$L_j(s):=\sum_{n\geq 1}a_{n^j} n^{-s}$.
Let ~$L(s,{\rm sym}^jg)$ be the symmetric cube and $4$th power
$L$-function attached to ~$g$ for ~$j=3,4$ respectively. In \cite[Lemma 3.3]{LAU-LU}, Lau and L\"u proved
the following. 
$$
L_j(s) = L(s,{\rm sym}^jg)\left(\prod_{1\leq i \leq j/2}L(2s,{\rm sym}^{2j-4i}g)^{-1}\right)H_j(s), ~~~~~~~~{\rm Re}(s) >1,
$$
where ~$H_j(s)$~ converges absolutely in the half-plane Re$(s)>1/3$.
Using the above result together with the properties of symmetric power $L$-function, 
we get the analytic continuation of the Dirichlet series $L_j(s), j=3,4$ for the region Re$(s)>1/2$. 
By using Phragmen-Lindel\"of
principle we get the polynomial growth for ~$L_j(s), j=2,3,4$~ in the region Re$(s)>1/2$.

\noindent By using the theory of symmetric power $L$-functions and their Rankin-Selberg convolution (for definitions and basic properties see \cite[\S 3]{LAU-LU}), 
we prove that for any ~$j=2,3,4$~
the Dirichlet series ~$\sum_{n\geq 1}a_{n^j}^2 n^{-s}$~ 
has a pole at $s=1$. In \cite[\S 7]{LY}, using the known facts about the $L$-functions $L(s,{\rm sym}^jg),j=2,3,4$ and their Rankin-Selberg convolution Li and Young
have studied the analytic properties of the series $\sum_{n\geq 1}a_{n^2}^4n^{-s}$ in detail. Following this we can write the series 
~$\sum_{n\geq 1}a_{n^j}^2 n^{-s}$~ as product of
the Rankin-Selberg $L$-function $L(s,{\rm sym}^jg\times {\rm sum}^jg)$ and another Dirichlet series $U_j(s)$, where the series $U_j(s)$ can be analytically continued
beyond the line Re$(s)=1$. 
Since $L(s,{\rm sym}^jg\times{\rm sym}^jg)$ has a pole at $s=1$, we see
from the above decomposition that the series
$$
\sum_{n\geq 1}a_{n^j}^2 n^{-s}
$$
has also a pole at $s=1$.
Also we have $\alpha+r = j(7/64+\epsilon)+1/2 < 1, j=2,3,4$.
Thus the subsequences of coefficients
$\{a_{n^j}\}_{n=1}^{\infty}, j=2,3,4$ satisfies all the hypothesis of Theorem \ref{signchange} and therefore they 
are oscillatory (infinitely many sign changes).
\end{enumerate}

\bigskip

\noindent {\bf Acknowledgement:} The second author was supported by a research fellowship from the Council of Scientific and Industrial Research (CSIR).


\bigskip

\begin{thebibliography}{100}
\bibitem{choie-kohnen1}
Y. Choie and W. Kohnen, {\em Mellin transforms attached to a certain automorphic integral},
J. Number Theory {\bf 132} (2012), 301--313.
\bibitem{DKS}
S. Das, W. Kohnen and J. Sengupta, {\em On a convolution series attached to a Siegel Hecke cusp form of
degree $2$}, Ramanujan J. {\bf 33} (2014), 367--378.
\bibitem{secondorder}
N. Diamantis, M. Knopp, G. Mason and C. O'Sullivan, {\em $L$-functions of second order cusp forms},
Ramanujan J. {\bf 12} (2006), 327--347.
\bibitem{hulse}
T. A. Hulse, E. M. Kiral, C. I. Kuan and L. Lim, {\em The sign of Fourier coefficients of half-integral 
weight cusp forms}, Int. J. Number Theory {\bf 8} (2012), 749--762.
\bibitem{iwaniec}
H. Iwaniec, Topics in classical automorphic forms, AMS, 1997.
\bibitem{KP}
J. Kaczorowski and A. Perelli, {\em General $\Omega$-theorems for coefficients of $L$-functions}, Proc. Amer. Math. Soc. (in press).
\bibitem{kim-sarnak}
H. Kim and P. Sarnak, {\em Appendix 2: Refined estimates towards the Ramanujan and Selberg 
conjectures}, J. Amer. Math. Soc. {\bf 16} (2003), 175--181.
\bibitem{KKP}
M. Knopp, W. Kohnen and W. Pribitkin, {\em On the signs of Fourier coefficients of cusp forms}, 
The Ramanujan Journal {\bf 7} (2003), 269--277.
\bibitem{K07} 
W. Kohnen, {\em Sign changes of Hecke eigenvalues of Siegel cusp forms of genus two},
Proc. Amer. Math. Soc. {\bf 135} (2007), 997--999.
\bibitem{kohnen}
W. Kohnen, {\em On the growth of the Petersson norms of Fourier-Jacobi coefficients of Siegel cusp forms},
Bull. Lond. Math. Soc. {\bf 43} (2011), 717--720.
\bibitem{kohnen-sengupta}
W. Kohnen and J. Sengupta, {\em The first negative Hecke eigenvalue of a Siegel cusp form of genus two},
Acta Arith. {\bf 129} (2007), 53--62.
\bibitem{LS} 
H. Lao, A. Sankaranarayanan, {\em The average behaviour of Fourier coefficients of 
cusp forms over sparse sequences}, Proc. Amer. Math. Soc. {\bf 8}(2009), 2557--2565.
\bibitem{LAU-LU} 
Y. K. Lau and G. L\"u, {\em Sums of Fourier coefficients of cusp forms}, Quart. J. Math. 
{\bf 62}(2011), 687--716.
\bibitem{LY} 
X. Li and M. P. Young, {\em Additive twists of Fourier coefficients of symmetric-square lifts}, J. Number Theory 
{\bf 132}(2012), 1626-1640.
\bibitem{LU} 
G. S. L\"u, {\em On an open problem of Sankaranarayanan}, Science China Math. 
{\bf 53}(2010), 1319--1324.
\bibitem{MSV}
J. Meher, K. D. Shankhadhar and G. K. Viswanadham, {\em A short note on sign changes}, 
Proc. Indian Acad. Sci. {\bf 123} (2013), 315--320.
\bibitem{pribitkin1}
W. Pribitkin, {\em On the sign changes of coefficients of general Dirichlet series}, 
Proc. Amer. Math. Soc. {\bf 136} (2008), 3089--3094.
\bibitem{pribitkin2}
W. Pribitkin, {\em On the oscillatory behavior of certain arithmetic functions associated with 
automorphic forms}, J. Number Theory {\bf 131} (2011), 2047--2060.
\bibitem{weissauer}
R. Weissauer, {\em Endoscopy for $GSp(4)$ and the cohomology of Siegel modular threefolds},
Springer Lecture notes in Mathematics {\bf 1968} (2009).
\end{thebibliography}
\end{document}